\documentclass[12pt,reqno]{amsart}
\UseRawInputEncoding

\usepackage[all]{xy}
\usepackage{amsthm,array,amssymb,amscd,amsfonts,latexsym}

\usepackage{enumerate}
\usepackage{mathrsfs}

\theoremstyle{plain}

\theoremstyle{definition}

\newtheorem{nothing*}[theorem]{}
\newtheorem{subnothing*}[sub]{}

\theoremstyle{remark}

\newcommand{\cc}{\raise .4pt \hbox{{$\scriptstyle{\bullet}$}}}

\begin{document}

\title[Picard group of connected affine algebraic group]{Picard group\\ of connected affine algebraic group}
\author[Vladimir L. Popov]{Vladimir L. Popov}
\address{Steklov Mathematical Institute,
Russian Academy of Sciences, Gub\-kina 8,
Moscow 119991, Russia}
\email{popovvl@mi-ras.ru}

\begin{abstract} We prove that the Picard group of a connected affine algebraic group $G$ is isomorphic
to the fundamental group of the derived subgroup of the reductive algebraic group $G/{\mathscr R}_u(G)$, where ${\mathscr R}_u(G)$ is the unipotent radical of $G$.
 \end{abstract}

\maketitle

All algebraic varieties considered below are defined over a basic algebraically closed field $k$. We follow the point of view on algebraic groups accepted in \cite{B1991} and use the following notation.

  If $S$ is a connected semisimple algebraic group, then
$\widehat S$ is its universal cover, and $\pi(S)$ is the kernel
  of the canonical isogeny $\widehat S\to S$.
  If $G$ is a connected affine algebraic group and $H$ is its closed subgroup, then
    $\varepsilon_{G, H}\colon {\rm Hom}_{\rm alg}(H, \mathbb G_m)\to {\rm Pic}(G/H)$ is the homomorphism that maps each character $\chi\colon H\to \mathbb G_m$
    to  the class of the one-dimensional homogeneous vector bundle over $G/H$ determined by $\chi$ (see \cite[Thm.\,4]{P1974}). If $\varphi\colon X\to Y$ is a morphism of
   smooth
irreducible algebraic varieties, then
$\varphi^*\colon {\rm Pic}(Y)\to {\rm Pic}(X)$ is the associated with $\varphi$
homomorphism of Picard groups (see\;\cite[Chap.\,III, \S1, Sect.\,2]{S2007}).
  Recall that the derived subgroup of a connected reductive algebraic group is connected and semisimple (see \cite[Sects. I.2.2 and II.14.2]{B1991}).

 The purpose of this note is to prove the following theorem.

\vskip 2mm
{\sc Theorem.} {\it Let $G$ be a connected affine algebraic group, let ${\mathscr R}_u(G)$ be its unipotent radical, let $\varrho\colon G\to G/{\mathscr R}_u(G)$ be the canonical homomor\-phism, let $S$ be the derived subgroup of the connected reductive group $G/{\mathscr R}_u(G)$, and let $\iota\colon S\hookrightarrow G/{\mathscr R}_u(G)$ be the identical embedding. Then the following canonical homomorphisms are isomorphisms:
\begin{equation*}
{\rm Pic}(G)\xleftarrow{\varrho^*}{\rm Pic}(G/{\mathscr R}_u(G))\xrightarrow{\iota^*}
{\rm Pic}(S)\xleftarrow{\varepsilon_{\widehat{S}, \pi(S)}}{\rm Hom}_{\rm alg}(\pi(S), \mathbb G_m).
\end{equation*}
}
\vskip -3mm

{\sc Corollary.} {\it The group ${\rm Pic}(G)$ is canonically isomorphic to the group ${\rm Hom}_{\rm alg}(\pi(S), \mathbb G_m)$ and is noncanonically isomorphic to the group $\pi(S )$.
}

\vskip 2mm

{\sc Example.} Let $G={\rm GL}_n$. Then the group ${\mathscr R}_u(G)$ is trivial, and the derived group of the group $G$ is the semisimple group ${\rm SL}_n$. The latter is simply connected, so the group $\pi({\rm SL}_n)$
is trivial. Therefore,
by Theorem, the group ${\rm Pic}(G)$ is trivial. This agrees with the fact that the group variety of the group ${\rm GL}_n$ is isomorphic to an open subset of $\mathbb A^{n^2}$.

\vskip 2mm

The following lemma is used in the proof of Theorem.

\vskip 2mm

{\sc Lemma.} {\it  Let $X$ be an irreducible smooth algebraic variety, let $U$ be a nonempty open subset of $\mathbb A^{d}$, and let $u_0$ be a point of $U$.
 Then for the morphisms
 \begin{equation}\label{new}
\begin{split}
\alpha\colon X\times U
&\to X,\quad (x, u)\mapsto x,\\
\beta\colon X
&\to X\times U,\quad x\mapsto (x, u_0),
\end{split}
\end{equation}
the homomorphisms $\alpha^*$ and $\beta^*$ are mutually inverse isomorphisms. }

\vskip 2mm

{\sc Proof of Lemma.} Consider  the morphisms
\begin{equation}\label{mo}
\begin{split}
\gamma\colon X\times \mathbb A^d
&\to X,\quad (x, a)\mapsto x,\\
\delta\colon X\times U
&\to X\times\mathbb A^d,\quad (x, u)\mapsto (x, u).
\end{split}
\end{equation}
It follows from \eqref{new} and
\eqref{mo} that the following equalities hold:
\begin{align}
\gamma\circ\delta\circ\beta
={\rm id}_X\qquad\mbox{and}\qquad
\alpha\circ\beta
={\rm id}_X.\label{two}
\end{align}
As is known, $\gamma^*$ is an isomorphism (see \cite[Chap.\,II, Prop.\,6.6 and its proof]{H1981}), and $\delta^*$ is a surjection (see \cite[Chap.\,II, Prop.\,6.5(a)]{H1981}). From this and the left equality in \eqref{two} it followes that $\beta^*$ is an isomorphism. In view of the right equality in \eqref{two}, this shows that $\alpha^*$ is also an isomorphism,
and  $\alpha^*$ and $\beta^*$ are mutually inverse.

\vskip 2mm

{\sc Proof of Theorem and Corollary.} Since the connected affine algebraic group ${\mathscr R}_u(G)$ is unipotent, it follows from \cite[Props.\,1, 2]{G1958} that

(a) the group variety of ${\mathscr R}_u(G)$ is isomorphic to an affine space,

(b) there is the commutative diagram
\begin{equation}\label{5}
    \xymatrix{
    G\ar[r]^(.23){\tau}\ar@/_1.4pc/[rr]^(.47){\varrho}&{\mathscr R}_u(G)\times \big(G/{\mathscr R}_u(G)\big)\ar[r]^(.62){\upsilon}&G/{\mathscr R}_u(G),
    }
    \end{equation}
    where $\tau$ is an isomorphism of group varieties (but, generally speaking, not of groups) and $\upsilon$ is the natural projection onto the second factor.

In view of Lemma, it follows from (a) and \eqref{5} that $\varrho^*$ is an isomorphism.

According to \cite[Thm.\,1]{P2022}, in the group $G/{\mathscr R}_u(G)$ there exists a torus $Z$ such that the mapping
$$\mu\colon S\times Z\to G/{\mathscr R}_u(G),\quad (s, z)\mapsto sz,$$
\noindent is an isomorphism of group varieties (but, in general, not of groups). Consider the commuta\-tive diagram
\begin{equation}\label{6}
\xymatrix{
S\ar[r]^(.4){\nu}\ar@/_1.3pc/[rr]^{\iota}&S\times Z\ar[r]^(.43){\mu}&G/{\mathscr R}_u(G)},
\end{equation}
in which
$\nu\colon S\to S\times Z$,\; $s\mapsto (s, e)$, where $e$ is the identity element.
Since the group variety of the torus $Z$ is isomorphic to an open subset of the affine space,
\eqref{6} and Lemma imply that $\iota^*$ is an isomorphism.

In view of the semisimplicity of the group
$\widehat S$ the group ${\rm Hom}_{\rm alg}(\widehat S, \mathbb G_m)$
is trivial, and since the group $\widehat S$ is simply connected,
the group ${\rm Pic}(\widehat S)$ is trivial (see
\cite[Prop.\,1]{P1974}).
According to \cite[Thm.\,4]{P1974}, it follows from this that
$\varepsilon_{\widehat{S}, \pi(S)}$ is an isomorphism. This completes the proof of Theorem.

The first part of Corollary follows directly from Theorem, while the second part follows from the fact that $\pi (S)$ is a finite abelian group.

\vskip 2mm

{\sc Remark.} The above Theorem
corrects Theorem 6 of \cite{P1974}. The latter asserts that the group ${\rm Pic}(G)$  is isomorphic to $\pi\big(G/{\mathscr R}(G)\big)$, where
${\mathscr R}(G)$ is the solvable radical of $G$.
If the group extension $1\to {\mathscr R}(G)\to G\to G/{\mathscr R}(G)\to 1$ splits,
then the group $G/{\mathscr R}(G)$ is isomorphic to the derived group of $G/{\mathscr R}_u(G)$, and so the formulated assertion is true in view of Theorem proved above.
But in general this is not the case, as Example above shows:
in it, the group $G/{\mathscr R}(G)$ is isomorphic to ${\rm PGL}_n$, and $\pi({\rm PGL}_n)$ is isomorphic
to the group of all $n$-th roots of $1$ in the field $k$. This latter group is nontrivial if $n$ is not a power of the characteristic of the field $k$ (however, the group ${\rm Pic}(G)$ is trivial for every $n$).

\vskip 2mm

I am grateful to Shuai Wang for bringing this example to my atten\-tion; this led to writing of this note. I am also indebted to S.\,O.\,Gor\-chin\-sky whose comments led to the above proof of Lemma and emphasis on the canonical nature of the construction
(the original proof of Lemma in preprint \cite{P2023}
was more geometric).

\end{document}